\documentclass[british,a4paper]{amsart}
\usepackage{babel}
\usepackage[T1]{fontenc}
\usepackage{amscd}
\usepackage{url}

\theoremstyle{plain}
\newtheorem{thm}{Theorem}[section]
\newtheorem{cor}[thm]{Corollary}
\newtheorem{prop}[thm]{Proposition}
\newtheorem{lem}[thm]{Lemma}

\theoremstyle{definition}
\newtheorem{df}{Definition}[section]

\theoremstyle{remark}
\newtheorem*{rmk}{Remark}
\newtheorem*{acks}{Acknowledgments}

\newtheorem*{claim}{Claim}

\newcommand{\ie}{\textit{i.e. }}

\newcommand{\p}{\mathbb{P}}

\newcommand{\N}{\mathcal{N}}
\newcommand{\Gr}{\mathbb{G}}

\newcommand{\OO}{\mathcal{O}}
\newcommand{\ud}{\mathrm{d}}

\newcommand{\abs}[1]{\lvert#1\rvert}

\newcommand{\Tt}{\mathcal{T}}
\newcommand{\Hh}{\mathcal{H}}

\DeclareMathOperator{\rk}{rk}

\DeclareMathOperator{\length}{length}

\DeclareMathOperator{\red}{red}

\DeclareMathOperator{\ndeg}{(\nu+1)-deg}

\DeclareMathOperator{\Hom}{Hom}

\DeclareMathOperator{\bl}{Bl}

\begin{document}

\title{Threefolds of $\p^5$ with one apparent quadruple point}

\author{Pietro De Poi}
\dedicatory{\textup{Dipartimento di Scienze Matematiche\\
Universit\`a degli Stud\^\i\ di Trieste\\
Via Valerio, 12/b\\
I-34127 Trieste\\ 
Italy\\
\url{depoi@mathsun1.univ.trieste.it}}} 
\thanks{The author has a Post-Doc grant at the Trieste University 
and is a member of GNSAGA-CNR} 


\begin{abstract}
In this article we classify all the smooth threefolds of $\p^5$ 
with an apparent quadruple point provided that the family of its 
$4$-secant lines is an irreducible (first order) congruence. This is 
sufficient to conclude the classification of all the smooth codimension 
two varieties of $\p^n$ with one apparent $(n-1)$-point and with 
irreducible family of $(n-1)$-secant lines. 
\end{abstract}

\maketitle

\bibliographystyle{amsalpha}

\section*{Introduction}
A congruence of lines in $\p^n$ is a family of lines of 
dimension $n-1$, and its order is the number of lines passing 
through a general point of $\p^n$. A codimension two subvariety of $\p^n$ is 
said to have $q$ apparent $(n-1)$-tuple points if its general projection from 
a point to a hyperplane has $q$ $(n-1)$-tuple points as singularities. 

In this work we prove that the degree of these 
varieties is bounded by $(n-1)^2$ when they have only one apparent 
$(n-1)$-tuple point (\ie $q=1$), and we observe that this implies that 
they cannot be complete intersections. 
By a result of A. Holme and M. Schneider, \cite{AM}, this implies also that 
in order to classify the smooth ones, we can stop up to dimension three. 

In this paper we also give a partial result towards this classification, \ie 
we restrict ourselves to the case in which the family of the 
$(n-1)$-secant lines is in fact an irreducible first order congruence. 

This article is structured as follows: 
after giving, in Section~\ref{sec:1}, the 
basic definitions, we obtain some general results about 
first order congruences in $\p^n$. 
In particular, after giving the central definition of 
fundamental $d$-loci, 
we show how to obtain the 
degree bound for the fundamental $(n-2)$-locus of a congruence. 

In Section~\ref{sec:6}, we give two general examples of congruences in $\p^n$: 
the first one is that of linear congruences, \ie the congruences which come 
out from general linear sections of the Grassmannians; the second one is 
that of congruences given by the $(n-1)$-secant lines of the varieties given 
by the degeneracy locus of a general map 
$\phi\in\Hom(\OO_{\p^n}^{\oplus (n-1)},\OO_{\p^n}^{\oplus n}(1))$. 
Of these two examples we calculate locally free resolutions of the ideal 
sheaves of their focal locus and  of the congruences themselves. 
We study these examples because they give us all the congruences of our 
classification but one (\ie case \eqref{case} of Theorem~\ref{thm:primo}). 

In Section~\ref{sec:2} we prove two multiple point formulae: the quadruple 
point formula for a smooth threefold of $\p^5$ and the formula which gives 
the number of $4$-secant lines to a smooth surface of $\p^4$ passing through 
a general point of the surface itself. 
Strangely enough, the quadruple point formula for the threefolds---which is 
actually an application of S. Kleiman's multiple point formulae of maps 
(see \cite{K1})---it seems to have been unknown, at least in modern 
times (see for example \cite{BSS}). With this formula, E. Mezzetti has been 
able to exclude the only 
degree $12$ smooth threefold of $\p^5$ for which the existence was uncertain 
(see \cite{Ed}). I must say that after the proof of this formula, I realized 
that the same was obtained by S. Kwak in \cite{Kw}, but with other 
methods, \ie through the monoidal construction. 

The irreducible congruences of order one of $\p^5$ 
which are given by the $4$-secant lines of a smooth 
threefold are classified in 
Section~\ref{sec:3}. 
By the results of Section~\ref{sec:6} and what we said above, 
we obtain the following complete list, where $d$ is the degree of the smooth 
codimension two variety $X$, $\pi$ its sectional genus and, if $\dim X=3$, 
$S$ is its general hyperplane section; finally, $H$ and $K$ are the hyperplane 
and canonical divisor classes, respectively (and $X_i$ is referred to the 
classification given in \cite{DP}): 

\begin{thm}\label{thm:primo} 
The smooth codimension two (irreducible) subvarieties of $\p^n$ for which the 
family of their $(n-1)$-secant lines is an irreducible first order congruence 
are 
\begin{enumerate}
\item for $n=2$ a point, and the congruence is a pencil of lines; 
\item for $n=3$ the twisted cubic, and the congruence is the Veronese 
surface and has bidegree $(1,3)$ (for more details, see 
Subsection~\textup{\ref{ssec:lin}}); 
\item for $n=4$ we have the following possibilities: 
\begin{enumerate} 
\item a (projected) Veronese surface, which is rational, with $d=4$, $\pi=0$; 
in this case we have a linear congruence, which has bidegree $(1,2)$ (see 
Subsection~\textup{\ref{ssec:lc}}); 
\item a Bordiga surface, which is rational, with $d=6$, $\pi=3$; 
the congruence is smooth and has bidegree $(1,8)$ (see 
Subsection~\textup{\ref{ssec:lin}}); 
\end{enumerate}
\item for $n=5$ we have the following possibilities: 
\begin{enumerate}
\item the Palatini scroll, which is rational with $d=7$, $\pi=4$, 
$\chi(\OO_S)=\chi(\OO_X)=1$ (case $X_6$); 
in this case we have a linear congruence, which has $3$-degree $(1,3,2)$ (see 
Subsection~\textup{\ref{ssec:lc}}); 
\item\label{case} a non rational scroll, $\p^1$-bundle over a minimal $K3$ 
surface of 
$\p^8$ via $\abs{K+H}$ (case $X_{11}$); $d=9$, $\pi=8$, 
$\chi(\OO_S)=\chi(\OO_X)=2$; the congruence, has $3$-degree $(1,7,13)$; 
\item a log-general type rational threefold, linked with a $(4,4)$ complete 
intersection to a Bordiga threefold (case $X_{15}$); 
$d=10$, $\pi=11$, $\chi(\OO_S)=5$, $\chi(\OO_X)=1$; 
the congruence is smooth (see 
Subsection~\textup{\ref{ssec:lin}}) and has $3$-degree $(1,15,20)$. 
\end{enumerate}
\end{enumerate}
\textit{Vice versa} the $(n-1)$-secant lines of 
any of the above varieties generate a first order congruence. 
\end{thm}

We conjecture that the first order congruences whose pure fundamental locus 
$F$ (see Definition \ref{df:pure} below) 
is a smooth variety of codimension two are the ones obtained in 
Theorem~\ref{thm:primo} (\ie the families of $(n-1)$-secant 
lines of the varieties of the theorem). 

\begin{acks}
I would like to thank Emilia Mezzetti, Ragni Piene, Kristian Ranestad, 
Fyodor~L.~Zak and Jon Eivind Vatne for interesting 
comments and remarks. 

I also thank the referee for the careful reading and important suggestions 
and remarks. 
\end{acks}

\section{Notations, definitions and general results}\label{sec:1}

We will work with schemes and varieties over the complex field $\mathbb C$. 
By \emph{variety} we mean a reduced and irreducible algebraic 
$\mathbb C$-scheme. More information about general results and references 
about families of lines, focal diagrams and congruences can be found in 
\cite{DP2} or \cite{PDP}. Besides, we refer to \cite{GH} for notations about 
Schubert cycles and to \cite{Fu} for the definitions and results of 
intersection theory. 
Here we recall that a \emph{congruence of lines} of $\p^n$ 
is a flat family $(\Lambda,B,p)$ of lines of $\p^n$ obtained by the 
desingularization of a subvariety $B'$ of dimension $n-1$ of the Grassmannian 
$\Gr(1,n)$ of lines of $\p^n$. $p$ is the restriction of the projection 
$p_1:B\times\p^n\rightarrow B$ to $\Lambda$, while we will denote the 
restriction of $p_2:B\times\p^n\rightarrow\p^n$ by $f$. 
$\Lambda_b:=p^{-1}(b)$, $(b\in B)$ will be an element of the family and 
$f(\Lambda_b)=:\Lambda(b)$ is a line of $\p^n$. 
We can summarise all these notations in the following two diagrams: 
the first one defines the family 
$$
\begin{CD}
\Lambda:=\psi^*(\Hh_{1,n}) @>\psi^*>> 
\Hh_{1,n}@>p_2>>\p^n\\
 @VpVV  @Vp_1VV \\
 B@>\psi >>B'\subset\Gr(1,n),
\end{CD}
$$
where $\Hh_{1,n}\subset\Gr(1,n)\times \p^n$ is the incidence variety 
and $\psi$ is the 
desingularization map, and the second one explains the notation for the 
elements of the family 
$$
\begin{CD}
\Lambda_b\subset\Lambda@>f:=\psi^* p_2>> 
\p^n\supset\Lambda(b):=f(\Lambda_b)\\
 @VpVV   \\
b\in B.
\end{CD}
$$
A point $y\in \p^n$ is called 
\emph{fundamental} if its fibre $f^{-1}(y)$ has dimension greater than the 
dimension of the general one. The \emph{fundamental locus} is the set of the 
fundamental points. 
The \emph{subscheme of the foci of the first order}, or, simply, 
the \emph{focal subscheme} 
$V\subset \Lambda$ is the scheme of the ramification points  of $f$. 
The \emph{locus of the first order foci}, or, simply, the \emph{focal locus}, 
$\Phi:=f(V)\subset \p^n$, is the set of the 
branch points of $f$. In this article, as we did in \cite{DP2}, 
we will endow this locus with the 
scheme structure given by considering it as the scheme-theoretic image of 
$V$ under $f$ (see, for example, \cite{H}). 

To a congruence is associated a \emph{sequence of degrees} or 
\emph{$(\nu+1)$-degree} 
$(a_0,\dotsc,a_\nu)$ if we write 
$$
[B]=\sum_{i=0}^{\nu}a_i\sigma_{n-1-i, i}
$$
---where we put $\nu:=\left [\frac{n-1}{2}\right ]$---as 
a linear combination of 
Schubert cycles of the Grassmannian (recall that 
$\sigma_{n-1-i, i}:=[\{\ell\in \Gr(1,n)\mid\ell\cap\p^i\neq\emptyset;\ 
\ell,\p^i\subset \p^{n-i}\}]$). In particular, the \emph{order} $a_0$ 
is the number of lines of $B$ passing through a general point of $\p^n$. 
The fundamental locus is contained in the focal locus and the two loci 
coincide in the case of a first order congruence, \ie through a focal 
point there will pass infinitely many lines of the congruence. 
An important result---independent of order and class---is the following 
(see also \cite{DP2}): 

\begin{prop}\label{prop:fofi}
On every line $\Lambda_b\subset \Lambda$ of the family, 
the focal subscheme $V$ 
either coincides with the whole  $\Lambda_b$---in which case $\Lambda(b)$ 
is called \emph{focal line}---or is a zero dimensional 
subscheme of $\Lambda_b$ of length $n-1$. Moreover, in the latter case, 
if $\Lambda$ is a first order congruence, $\Phi\cap \Lambda_b$ has length  
$n-1$. 
\end{prop}

\begin{proof}
Let 
$\lambda:\Tt_{(B\times \p^n / \p^n)\vert_\Lambda}\rightarrow 
\N_{\Lambda/ B\times \p^n}$ be the \emph{global characteristic map} 
for the family $\Lambda$ (see \cite{CS}). From the 
\emph{focal diagram} (diagram (3) of  \cite{CS}) one gets that 
the subscheme of the foci of the first order $V$ is the degeneracy locus 
of $\lambda$. If we restrict the map $\lambda$ to a fibre 
$\Lambda_b\cong\Lambda(b)$, 
we obtain the \emph{characteristic map of the family relative to $b$}: 
$$
\begin{CD}
\lambda(b):T_{B,b}\otimes\OO_{\Lambda(b)} @>>> \N_{\Lambda(b)/ \p^n}\\
 @V\cong VV  @V\cong VV \\
 \OO^{n-1}_{\Lambda(b)}@>>>\OO_{\Lambda(b)}(1)^{n-1}.
\end{CD}
$$
From the preceding isomorphisms, the map 
$\lambda(b)$ can be seen as an $(n-1)\times(n-1)$-matrix with linear 
entries on $\Lambda(b)$; so the focal locus on $\Lambda(b)$ is 
given by the vanishing of the determinant of this matrix, and our claim 
follows. 

Concerning the first order congruences, we observe that a fundamental point 
$P$ is a focal point for every line $\Lambda(b)$ which contains it, since 
the characteristic map relative to it, $\lambda(b)$, drops rank in $P$. 
\end{proof}

\begin{rmk}
If $\Lambda$ is a congruence of order different from one and the intersection 
$F\cap\Lambda(b)$ is proper, this can have length greater than $n-1$. 
This is due to the fact that $\Lambda(b)$ can contain points of $\Phi$ that 
are actually focal for other lines, \ie points that are images of points 
contained in some $V\cap\Lambda_{b'}$, with $b\neq b'$. 
\end{rmk}


From now on we will consider only first order congruences, and since in this 
case the focal and the fundamental loci coincide, we will talk only of 
fundamental loci (with the scheme structure given above for the focal locus). 
A central definition, introduced first in \cite{DP2}, is that of 
\emph{fundamental $d$-locus}, \ie the subscheme of 
the fundamental locus of pure dimension $d$, with $0\le d\le n-2$, 
which is met by the general line of the congruence. 
Let us see how these schemes are constructed: 
the closed set 
$$
S_d:=\{(\Lambda(b),P)\in\Lambda\mid \rk(\ud f_{(\Lambda(b),P)})\le d\}
$$
has a natural subscheme structure, 
which is defined by a Fitting ideal, \ie 
the ideal generated by the $(d+1)$-minors of $\ud f$ (or, by the 
Fitting lemma, see \cite{E}, by the $d$-minors of $\lambda$), 
see \cite{SK}; in particular, $S_{n-1}=V$. 
Let us define 
$$
D_{d+1}:=\overline{S_{d+1}\setminus S_d}
$$ 
with the scheme structure induced by $S_{d+1}$. Finally, we consider the 
scheme-theoretic image $\Phi_d$ of $D_{d+1}$ in $\p^n$ under $f$. 
The component of $\Phi_d$ of pure dimension $d$ (with the scheme structure 
induced by $\Phi_d$) which is met by the general 
line of the congruence is the fundamental $d$-locus. 

These subschemes of the fundamental locus are particularly important, 
since a first order congruence can be characterised as a component of 
the set of lines which meet 
the fundamental $d$-loci a certain number of times, see the 
Classification Theorem~\textup{3.2} of \cite{DP2}. 

\begin{df}\label{df:pure} The union  of the fundamental $d$-loci of  $F$ 
is called \emph{pure fundamental locus},
or, in what follows, simply  \emph{fundamental locus} and it is denoted 
by $F$. 
\end{df}


After this, we give the following 

\begin{thm}\label{thm:dgb}
If $\Lambda$ is a first order congruence such that the pure fundamental locus 
$F$ is irreducible and coincides with the fundamental $(n-2)$-locus, then 
\begin{equation*}
\frac{n-1}{k}<m<(n-1)^2, 
\end{equation*}
where $m:=\deg (F)_{\red}$ and $k$ is the geometric multiplicity 
$(F)_{\red}$ in $F$. 
\end{thm}

\begin{proof}
First of all, we have that $n-1<km$ by degree reasons, since the congruence 
is given by lines which intersect $F$ in a zero dimensional scheme of 
length  $n-1$.

To prove the other bound, we need of course more work. 
Let $B$ be our congruence, which has sequence of degrees 
$(1,a_1,\dotsc,a_\nu)$. If $\Pi$ is a (fixed) general 
$(n-2)$-plane, we denote by $V_\Pi$ the scroll 
given by the lines of the congruence that 
meet $\Pi$. Then by the Schubert calculus one can show that 
$V_\Pi$ is a hypersurface of $\p^n$ of degree $1+a_1$. 

Moreover, if $\ell$ is a line of $B$ 
not contained in $V_\Pi$ and $P$ is a point of $V_\Pi\cap \ell$, 
then $P$ is a focus for $B$, since at least two lines of the 
congruence pass through it. 

Then, if $\Pi'$ is another general 
$(n-2)$-plane of $\p^n$, the complete intersection 
of the hypersurfaces 
$V_\Pi$ and  $V_{\Pi'}$ 
is a (reducible) $(n-2)$-dimensional scheme $\Gamma$ that contains 
the focal locus $F$ and 
the $(n-2)$-dimensional scroll 
$\Sigma$ given by the lines of the congruence 
meeting $\Pi$ and $\Pi'$, 
which has degree $1+2a_1+a_2$ (actually, $1+a_1$ in the case of $\p^3$; 
since this case can be treated 
analogously, we will suppose from now on that $n>3$). 

In fact, if a point $P$ of $V_\Pi\cap V_{\Pi'}$ 
does not belong to the scroll 
$\Sigma$, it belongs to the fundamental locus. 
Indeed in this case 
$$
P\in \ell\cap \ell', \ \text{where}\ \ell\in G_\Pi,\ \ell\in G_{\Pi'},
\ \text{and}\ \ell\neq \ell' 
$$
---where $G_\Pi$ and $G_{\Pi'}$ 
denote the subvarieties of the Grassmannian 
corresponding to the two scrolls 
$V_\Pi$ and  $V_{\Pi'}$. Since $\Lambda$ is a first order congruence and $P$ 
belongs to two of the lines of $\Lambda$, it belongs to infinitely many ones. 

Reciprocally, if $P\in F$ is a general focal point, 
the set of the lines of $B$ through $P$, 
$\chi_P$, is a cone of dimension (at least) two, so its intersection with 
$\Pi$ and $\Pi'$ is not empty and therefore $P\in V_\Pi\cap V_{\Pi'}$. 

The degree of the scroll follows from the Schubert 
calculus. 

\begin{claim}
The following formulae hold: 
\begin{align}
(n-1)h& =  1+a_1,\label{eq:agen}\\
(1+a_1)^2&
\ge h^2m+1+2a_1+a_2, 
\label{eq:bgen}
\end{align}
where we denoted by $h$ the algebraic multiplicity 
of $(F)_{\red}$ on $V_\Pi$. 
\end{claim} 

Let us start proving relation~\eqref{eq:agen}. 
If we take a line $\Lambda(b)$ of the 
congruence not contained 
in $F\cap V_\Pi$, then, intersecting $\Lambda(b)$ with $V_\Pi$, we obtain a 
zero dimensional scheme of length $1+a_1$, since this is 
the degree of $V_\Pi$. 
$V_\Pi\cap \Lambda(b)$, set-theoretically, coincides with 
$(F)_{\red} \cap\Lambda(b)$, and its 
intersection multiplicity is $h$, 
so the relation~\eqref{eq:agen} is proved. 

We recall that 
the degree of $\Gamma$ is $(1+a_1)^2$, 
and it contains $F$ and the scroll $\Sigma$. 
Actually, the fundamental $(n-2)$-locus has geometric 
multiplicity in $\Gamma$ equal to $h^2$. The proof of this fact is 
the following: 
the intersection multiplicity 
$i((F)_{\red},V_\Pi\cdot V_{\Pi'},\p^n)$ of $F$ in 
$V_\Pi\cdot V_{\Pi'}$ is equal to the geometric 
multiplicity of  $(F)_{\red}$ in $\Gamma$, 
but $i((F)_{\red},V_\Pi\cdot V_{\Pi'},\p^n)=h^2$. 

Finally, as we seen, the scroll $\Sigma$ has degree 
$1+2a_1+a_2$, so we get formula~\eqref{eq:bgen}.

Now, if we substitute formula~\eqref{eq:agen} 
in formula \eqref{eq:bgen}, we obtain 
\begin{equation}
(n-1)^2h^2 -h^2m-1-2a_1-a_2\ge 0,\label{eq:a2gen}
\end{equation}
and since $-1-2a_1-a_2 < 0$, we deduce $m<(n-1)^2$. 
\end{proof}

A fundamental consequence of the preceding theorem is the following: 

\begin{thm}\label{thm:HC}
If we have a first order congruence of $\p^n$ 
such that the fundamental locus $F$ 
satisfies the hypothesis of the preceding theorem, then $F$ cannot 
be a complete intersection. If moreover $F$ is smooth, then $n\le 5$. 
\end{thm}

\begin{proof}
In fact, by the preceding theorem $\deg(F)<(n-1)^2$, therefore if it were a 
complete intersection, it would be contained in a hypersurface $V$ 
of degree less than $n-1$, and so every $(n-1)$-secant line of $F$ would be 
contained in $V$. 

If $F$ is smooth, since by \cite{AM} we know that Hartshorne's conjecture is 
true in codimension two up to degree $(n-1)(n+5)$, then $\dim(F)\le 3$. 
\end{proof}


\section{General examples of first order congruences}
\label{sec:6}

We give now two examples of first order congruences 
of lines of $\p^n$. Actually, these examples give us all the congruences 
of Theorem \ref{thm:primo} but case \ref{case}.

\subsection{Linear sections of $\Gr(1,n)$}
\label{ssec:lc}
First of all, we will analyse the congruences that come out 
from linear sections of the Grassmannian $\Gr(1,n)$, 
\ie we will consider the so called, classically, 
\emph{linear congruences}.

We recall that the Schubert cycle that corresponds to a hyperplane section 
of (the projective embedding of) the Grassmannian is $\sigma_1$, so the 
following technical lemma gives us the formula for the general intersection 
of these special Schubert cycles: 

\begin{lem}
If $\ell\le n-1$ and we set $k:=\left[\frac{\ell}{2}\right]$, 
the following formula holds: 
\begin{equation}\label{eq:cie}
\sigma_1^\ell=\sum_{i=0}^{k}\left(\binom{\ell-1}{i}-\binom{\ell-1}{i-2}\right)
\sigma_{\ell-i,i}
\end{equation}
---with the convention that $\binom{\ell}{h}=0$ if $h<0$. 
\end{lem}

\begin{proof}
Let us prove the lemma by induction: for $\ell=1$ it is obvious. 
Let us suppose it is true for $\ell-1$; then, by inductive hypothesis 
\begin{equation*}
\sigma_1^{\ell-1}=\sum_{i=0}^{k'}
\left(\binom{\ell-2}{i}-\binom{\ell-2}{i-2}\right)\sigma_{\ell-1-i,i},
\end{equation*}
where $k':=\left[\frac{\ell-1}{2}\right]$. By Pieri's~formula, we have 
\begin{equation*}
\sigma_{\ell-1-i,i}\cdot \sigma_1 =
\begin{cases}
\sigma_{\ell-i,i}& \text{if $\ell-1-i=i$}\\ 
\sigma_{\ell-i,i}+\sigma_{\ell-1-i,i+1} & \text{otherwise,}
   \end{cases}
\end{equation*}
\ie if $\ell-1\neq 2k'+1$ and $i\neq k'$, we obtain 
$$
\binom{\ell-2}{i}-\binom{\ell-2}{i-2}+\binom{\ell-2}{i-1}-\binom{\ell-2}{i-3}=
\binom{\ell-1}{i}-\binom{\ell-1}{i-2}
$$
while if $i=k'=\frac{\ell-2}{2}$, 
$$
\binom{2k'}{k'}-\binom{2k'}{k'-2}=
\binom{2k'+1}{k'+1}-\binom{2k'+1}{k'-1}, 
$$
which follows from the formula 
\begin{equation}\label{eq:new}
\binom{\ell-1}{i}-\binom{\ell-1}{i-2}
=\binom{\ell}{i}\cdot\frac{\ell-2i+1}{\ell-i+1}. 
\end{equation}
\end{proof}

\begin{thm}
If $\Lambda$ is a congruence with sequence of degrees $(a_0,\dotsc,a_\nu)$ 
then $B$, as a subvariety of the Pl\"ucker embedding of the Grassmannian, has 
degree 
\begin{equation*}
\deg(B)=\sum_{i=0}^\nu a_i\left(\binom{n}{i}\cdot\frac{n-2i+1}{n-i+1}\right).
\end{equation*}
\end{thm}

\begin{proof}
It is a corollary of formulae~\eqref{eq:cie} and \eqref{eq:new}. 
\end{proof}

\begin{cor}\label{cor:ls}
An $(n-1)$-linear section $B$ of the Grassmannian of lines of $\p^n$ 
generates a first order congruence $\Lambda$ with sequence of degrees 
$$
\ndeg(\Lambda)
=\left(1,\dotsc,\left(\binom{n-2}{i}-\binom{n-2}{i-2}\right),\dotsc,
\left(\binom{n-2}{\nu}-\binom{n-2}{\nu-2}\right)\right);
$$
in particular, as a subvariety of the Pl\"ucker embedding of 
the Grassmannian, this is a smooth congruence of degree 
$$
\deg(B)=\sum_{j=0}^{\nu}\left(\binom{n-2}{j}-\binom{n-2}{j-2}\right)^2=
\frac{1}{n-1}\binom{2n-2}{n}. 
$$
\end{cor}

Obviously the degree of a linear congruence is the degree of the Grassmannian 
(which can be found in \cite{Fu}, Example 14.7.11). 
This corollary gives us a first non-trivial example of a first 
order congruence. Some general results about fundamental varieties of these 
congruences are given in \cite{BM}; in particular, it is proven that 
the focal locus is the degeneracy locus $F$ of a general morphism 
\begin{equation}\label{mf}
\phi:\OO_{\p^n}^{\oplus (n-1)}\rightarrow \Omega_{\p^n}(2)
\end{equation}
of (coherent) sheaves on $\p^n$ and that $F$ is smooth if 
$\dim(F)\le 3$. 

The idea of the proof is the following: first we have the isomorphism: 
\begin{equation}\label{com}
H^0(\Omega_{\p^n}(2))\cong(\wedge^2 V)^*,
\end{equation}
where $V:=H^0(\OO_{\p^n}(1))^*$; 
from this we can interpret a global section of $\Omega_{\p^n}(2)$ as 
a bilinear alternating form on $V$, or as a skew-symmetric matrix 
of type $(n+1)\times(n+1)$ with entries in the base field. 
Then, the general morphism $\phi$ defined in \eqref{mf} is assigned by giving 
$n-1$ general skew-symmetric matrices, $A_1,\dotsc,A_{n-1}$ and the 
corresponding 
degeneracy locus $F$ in $\p^n$ is the set of points $P$ such that 
\begin{equation}\label{sks}
\sum_{i=1}^{n-1}\lambda_i A_i[P]=0
\end{equation}
for some $(\lambda_1,\dotsc,\lambda_{n-1})\neq (0,\dotsc,0)$, and $[P]$ 
denotes the column matrix of the coordinates of $P$. 
Finally, recalling that the Pl\"ucker embedding of the Grassmannian is: 
\begin{equation*}
\psi:\Gr(1,n)\hookrightarrow \p(\wedge^2V)
\end{equation*}
we can conclude easily.

An improvement of a result of \cite{BM} is the following: 
\begin{prop}\label{prop:flc}
If $F$ is the focal locus of a general linear congruence of $\p^n$, then 
\begin{enumerate}
\item if $n$ is even, equation \eqref{sks} has at least one solution and 
$F$ is a rational variety; 
\item if $n$ is odd, 
the vanishing of the Pfaffian of the matrix of \eqref{sks} defines a 
hypersurface $Z$ of degree $(n+1)/2$ contained in a $\p^{n-2}$ 
in which $\lambda_1,\dotsc,\lambda_{n-1}$ are the coordinates. Furthermore, if 
$\phi$ is general, for a fixed point $[\lambda]\in Z$, we find a line of 
solutions of equation \eqref{sks} in $F$, so that $F$ is a scroll 
over (an open set of) $Z$. 
\end{enumerate} 
Besides, 
\begin{equation*}
\deg(F)=\frac{n^2-3n+4}{2}.
\end{equation*}
\end{prop}

\begin{proof}
If we dualise the Eagon-Northcott complex applied to \eqref{mf}, we have the 
following exact sequence 
\begin{equation*}
0\rightarrow\OO_{\p^n}(1-n)\rightarrow\Tt_{\p^n}(-2)
\xrightarrow{{}^t\phi}\OO_{\p^n}^{\oplus (n-1)}
\rightarrow\omega^\circ_F(2)\rightarrow 0,
\end{equation*}
where $\omega^\circ_F$ is the dualising sheaf of $F$. 
Hunting in the sequence, we get
$$
H^0(\p^n,\omega^\circ_F(2))\cong H^0(\p^n,\OO_{\p^n}^{\oplus (n-1)}),
$$
and so the complete linear system $\abs{K_F+2H}$---where $K_F$ and $H$ are the 
canonical and hyperplane divisors, respectively---is base point free, 
from which we obtain that 
the map associated to it  
$$
\varphi_{\abs{K_F+2H}}:F\rightarrow \p^{n-2} 
$$
is well defined. The fibre of  $\varphi_{\abs{K_F+2H}}$---if we fix projective 
coordinates $\lambda_1,\dotsc\lambda_{n-1}$ on $\p^{n-2}$ and 
$x_0\dotsc x_n$ on $\p^n$---is given, in view of \eqref{sks},
by the solutions of the following homogeneous linear system of $n+1$ 
equations in the $n+1$ 
indeterminates: 
$x_0\dotsc x_n$:
\begin{equation}
\sum_{\substack{i=1,\dotsc,n-1\\
k=0,\dotsc,n}}\lambda_i (A_i)_{j,k}x_k=0, \quad j=0,\dotsc,n,
\end{equation}
where $\lambda_1,\dotsc\lambda_{n-1}$ are fixed. Clearly the 
$(n+1)\times(n+1)$-matrix $A:=\sum \lambda_i (A_i)$ associated to this system 
is antisymmetric; 
therefore, if $n$ is even, its determinant is zero and has only one or 
infinitely many (projective) solutions; so, for dimensional reasons, we get 
that $\varphi_{\abs{K_F+2H}}$ is birational and $F$ is rational.

If instead $n$ is odd, $\det A$ is, in general, not zero, and so its Pfaffian 
defines a hypersurface $Z$ of degree $(n+1)/2$ and the general fibre of 
$\varphi_{\abs{K_F+2H}}$ is a $\p^1$. We observe that the Eagon-Northcott 
gives in fact a locally free resolution of the ideal sheaf of $F$: 
\begin{equation*}
0\rightarrow\OO_{\p^n}^{\oplus (n-1)}(1-n)
\xrightarrow{\phi(1-n)}\Omega_{\p^n}(3-n)
\rightarrow\OO_{\p^n}\rightarrow\OO_F\rightarrow 0.
\end{equation*}
From this we get the Hilbert polynomial of $F$, and 
in particular, after some computations, the degree. 
\end{proof}

\begin{rmk}
In low dimension and with a general section, we have that (see \cite{BM} for 
details) if $n=3$, $F$ is the union of two skew lines, if $n=4$, $F$ is a 
smooth projected Veronese surface and if $n=5$, $F$ is a (rational) threefold 
of degree seven, which is a scroll over a cubic surface in $\p^3$. It is also 
known as Palatini scroll (see \cite{O}).
\end{rmk}

\begin{rmk}
Since a linear congruence $B$ is given by a linear section of the 
Grassmannian $\Gr(1,n)$,  
a resolution of its ideal sheaf is given by the Koszul complex 
\begin{multline*}
0\rightarrow \OO_{\Gr(1,n)}(-n+1)
\rightarrow (\wedge^2((\OO_{\Gr(1,n)}(-1))^{\oplus (n-1)}))\rightarrow\dotsb\\
\dotsb\rightarrow((\OO_{\Gr(1,n)}(-1))^{\oplus (n-1)})
\rightarrow\OO_{\Gr(1,n)}\rightarrow\OO_B\rightarrow 0. 
\end{multline*}
Instead, the minimal free resolution of the Pl\"ucker embedding of $B$ can 
be obtained by the minimal free resolution of $\Gr(1,n)$---which can be 
found in \cite{JPW}---via the mapping cone. Since the computations are 
rather complicate, we will not give them. 
\end{rmk}


\subsection{Matrices of type $(n-1)\times n$ with linear entries}
\label{ssec:lin}

Let us consider a general morphism 
$\phi\in\Hom(\OO_{\p^n}^{\oplus (n-1)},\OO_{\p^n}^{\oplus n)}(1))$, whose 
minors vanish in the expected codimension two. In this case, 
$F:=V(\phi)$---the degeneracy locus of $\phi$---is a locally Cohen-Macaulay 
subscheme, the Eagon-Northcott complex is exact (see \cite{BE}) and gives a 
free resolution of our ideal sheaf:
\begin{equation}\label{eq:sh}
0\rightarrow \OO_{\p^n}^{\oplus (n-1)}(-n)\xrightarrow{\phi(-n)}
\OO_{\p^n}^{\oplus n}(1-n)\rightarrow \OO_{\p^n}\rightarrow \OO_F
\rightarrow 0. 
\end{equation}
Then---for example---from the Hilbert polynomial we get
\begin{align}
\deg(F)&=\binom{n}{2}\label{eq:deg}\\ 
\pi(F)&=1+\frac{2n-7}{3}\binom{n}{2}\label{eq:sg}
\end{align}
where $\pi(F)$ is the sectional genus of $F$. It is easy to prove that:
\begin{prop}\label{prop:pid}
$F$ is rational, and if $n\le 5$ it is smooth. Besides, the adjunction map 
$\varphi_{\abs{K_F+H}}$ exhibits $F$ as the blow-up of $\p^{n-2}$ in a scheme 
$Z$ of degree $\binom{n+1}{2}$ and sectional genus $\frac{n}{6}(2n-5)(n+1)-1$.
 In particular, if $n=4$, $F$ is a rational sextic which is the blow-up of the 
plane in $10$ points, \ie a \emph{Bordiga surface}. 
\end{prop}

\begin{proof}
The smoothness is a consequence of Bertini type theorems. For proving that $F$ 
is rational, we can apply the standard argument used in the proof of 
Proposition~\ref{prop:flc}: we consider the dual of 
the Eagon-Northcott complex applied to $\phi$
\begin{equation*}
0\rightarrow\OO_{\p^n}(-n)\rightarrow\OO_{\p^n}^{\oplus n}(-1)
\xrightarrow{{}^t\phi}\OO_{\p^n}^{\oplus (n-1)}
\rightarrow\omega^\circ_F(1)\rightarrow 0,
\end{equation*}
where $\omega^\circ_F$ is the dualising sheaf of $F$. 
Chasing in the sequence, we get
$$
H^0(\p^n,\omega^\circ_F(1))\cong H^0(\p^n,\OO_{\p^n}^{\oplus (n-1)}),
$$
and the complete linear system $\abs{K_F+H}$---where $K_F$ and $H$ are the 
canonical and hyperplane divisors, respectively---is base point free, 
from which we obtain that the map associated to it 
$$
\varphi_{\abs{K_F+H}}:F\rightarrow \p^{n-2} 
$$
is well defined. The fibre of $\varphi_{\abs{K_F+H}}$---if we fix projective 
coordinates $y_0,\dotsc y_{n-2}$ on $\p^{n-2}$ and 
$x_0\dotsc x_n$ on $\p^n$---is 
given by the solutions of the following homogeneous linear system of $n$ 
equations in the $n+1$ indeterminates $x_0\dotsc x_n$:
\begin{equation}
\sum_{i=0}^n \Phi_{j,i}(y_0:\dotsc:y_{n-2})x_i=0, \quad j=1,\dotsc,n,
\end{equation}
where $y_0,\dotsc y_{n-2}$ are fixed and $\Phi_{j,i}(y_0:\dotsc:y_{n-2})$ 
are linear forms on $\p^{n-2}$. Clearly the 
$(n+1)\times n$-matrix associated 
to this system 
has maximal rank, since the forms $\Phi_{j,i}(y_0:\dotsc:y_{n-2})$ 
are linearly independent if $\phi$ is general, and 
then it has only one (projective) solution; so, $\varphi_{\abs{K_F+H}}$ is 
birational and $F$ is rational. Besides, $\varphi_{\abs{K_F+H}}$ has (at 
least) one dimensional fibres on the degeneracy locus $Z$ of the map 
$\Phi:\OO_{\p^{n-2}}^{\oplus n}\rightarrow \OO_{\p^{n-2}}^{\oplus (n+1)}(1)$, 
\ie $\varphi_{\abs{K_F+H}}$ gives $F$ as the blow-up of $\p^{n-2}$ in $Z$. As 
usual, Eagon-Northcott gives a free resolution of $Z$:
\begin{equation*}
0\rightarrow\OO_{\p^{n-2}}^{\oplus n}(-n-1)\rightarrow
\OO_{\p^{n-2}}^{\oplus (n+1)}(-n)\rightarrow\OO_{\p^{n-2}}
\rightarrow\OO_Z\rightarrow 0,
\end{equation*}
from which we obtain the degree and the sectional genus of $Z$. 
\end{proof}

With this,  we prove that 

\begin{thm}\label{thm:bo1}
The $(n-1)$-secant lines of the variety $F$ defined as above form a first 
order congruence of lines $B$ of $\p^n$. The congruence $B$ is smooth for 
general $\phi$. 
\end{thm}

\begin{proof}
Let 
$$
A:=\begin{pmatrix}
a_{11}&\dotsb &a_{1 (n-1)}\\
\vdots&&\\
a_{n1}&\ &a_{n (n-1)}
\end{pmatrix}
$$
be the matrix of linear entries on $\p^n$ associated to the map $\phi$. 
Moreover, let us denote by $A_i$, $i=1,\dotsc, n$ the rows of $A$ and 
by $a_{\hat{i}}$ the minor of order $n-1$ of $A$ obtained cutting out the 
row $A_i$.  
We want to prove that any $(n-1)$-secant line of $F$ 
is given by the vanishing of the entries of 
$A_\lambda:=\sum_{i=1}^{n}\ell_{\lambda,i} A_i$, \ie by a linear 
combination of the rows of $A$: 
$\ell_\lambda:=V(\sum_{i=1}^{n}\ell_{\lambda,i} a_{ij})$, $j=1,\dotsc, n-1$. 
In fact, $\ell_\lambda$ is an $(n-1)$-secant line of $F$: after changing 
coordinates, we can assume that $\ell_\lambda$ is defined by the zeros of 
the last row. The intersection 
$V(a_{\hat{n}})\cap \ell_\lambda$ is given by $n-1$ points 
$P_1,\dotsc,P_{n-1}$ and by definition, the matrix $A$ drops rank at each of 
the $P_i$'s. 

\textit{Vice versa}, if $\ell$ is an $(n-1)$-secant line, to prove that its 
equations are a linear combination of rows it is sufficient to show that 
$\ell$ is contained in at least one (linear combination) of the hypersurfaces 
$V(a_{\hat{i}})$. If this were not so, we would have that 
$V(a_{\hat{i}})\cap\ell=F\cap\ell =\{P_1,\dotsc,P_{n-1}\}$, 
$\forall i$, since $V(a_{\hat{i}})\cap\ell\supset F\cap\ell$ and both 
sets contain $n-1$ points. But if we take a pencil of these hypersurfaces, 
we can find a linear combination of these containing $\ell$.  

By this, it is easy to see that $B$ is a first order congruence: 
if $P\in \p^n$ is a general point, we have to solve  
$A(P)\cdot \ ^t \lambda$ with $\lambda\in\p^{n-1}$, \ie a homogeneous 
linear system of $n-1$ equations in $n$ indeterminates.

Finally, for proving that $B$ is smooth, from the above geometric description 
of the $(n-1)$-secant lines as linear combination of the rows of $A$, 
we see that the map $\phi$ gives rise 
to a map 
\begin{equation*}
\varphi:\OO^{\oplus n}_{\Gr(1,n)}\rightarrow (\mathcal{S}^*)^{\oplus (n-1)}.
\end{equation*}
---where $\varphi$ is obtained by considering the dual of $\phi$ twisted by 
one and then pulled back to the incidence variety and finally pushed it 
forward to $\Gr(1,n)$---and $B$ is precisely the degeneracy locus of 
$\varphi$. Then we conclude as usual by Bertini type theorems. 
\end{proof}

\begin{rmk} 
If we apply the Eagon-Northcott complex to $\varphi$, we obtain a 
resolution of $B$:  
\begin{multline}\label{lira}
0\rightarrow\OO^{\oplus \binom{n}{2}}_{\Gr(1,n)}
\rightarrow(\mathcal{S}^*)^{\oplus (n-1)\binom{n}{3}}\rightarrow\\
\rightarrow (\wedge^2((\mathcal{S}^*)^{\oplus (n-1)}))^{\oplus \binom{n}{4}}
\rightarrow\dotsb\\
\dotsb\rightarrow\wedge^{n-2}((\mathcal{S}^*)^{\oplus (n-1)})
\rightarrow\OO_{\Gr(1,n)}(n-1)\rightarrow\OO_B(n-1)\rightarrow 0.
\end{multline}
\end{rmk}


\section{Two multiple point formulae}\label{sec:2}

We will prove now the quadruple point formula for a smooth 
threefold $X$ of $\p^5$ and the formula which gives the number of $4$-secant 
lines to a surface $S$ of $\p^4$ passing through a general point $P\in S$.  

In this and in the next section, 
we will denote by $S$ a general hyperplane section 
of the threefold $X$, by $H$ and $K$ the hyperplane and the canonical divisors 
of $X$. $d:=\deg(X)=H^3$ is the degree of $X$, 
$\pi=\frac{1}{2}H^2(K+2H)+1$ its sectional genus, 
while $\chi(\OO_X)$ and $\chi(\OO_S)$ are, respectively, 
the Euler-Poincar\'e characteristic of $X$ and $S$. 
We recall that for a smooth threefold of $\p^5$ two double point formulae 
hold---\ie one for $X$ and one $S$---which can be written as 
\begin{align*} 
K^3&=-5d^2+d(2\pi+25)+24(\pi-1)-36\chi(\OO_X)-24\chi(\OO_S),\\
H\cdot K^2&=\frac{1}{2}d(d+1)-9(\pi-1)+6\chi(\OO_X), 
\end{align*}
so the basic invariants of $X$ are $d,\pi,\chi(\OO_X)$ and $\chi(\OO_S)$ 
(see for example \cite{DP}).

We start with the quadruple point formula. We refer to \cite{K1} for the 
definitions and results used in the proof. 

\begin{prop}
Let $X$ be a smooth threefold of $\p^5$ and $P\in(\p^5\setminus X)$ is a point 
through which there is a finite number of $4$-secant lines of $X$.  
If we set $q(X):=\length(\Psi\cap W_P)$ where 
$\Psi, W_P\subset\Gr(1,n)\subset \p^{14}$ 
are, respectively, the family of $4$-secant lines of $X$ and the Schubert 
variety of the lines through $P$ (embedded via the Pl\"ucker embedding), 
then the following formula holds:
\begin{equation}
\begin{split}
q(X)& =\frac{1}{24}d^4-\frac{1}{4}d^3+\frac{1}{2}d^2(\frac{11}{12}-\pi)
+d(\frac{5}{2}\pi+2\chi(\OO_S)-\frac{9}{4})+\\
&\quad +\frac{1}{2}\pi^2-\frac{7}{2}\pi+6\chi(\OO_X)-9\chi(\OO_S)+3,  
\label{eq:4tuple}
\end{split}
\end{equation}
In particular, for the general point of $\p^5$, there pass  
$q(X)$ $4$-secant lines, \ie $\Psi$---if it is not empty---is a (possibly 
reducible) congruence of order $a_0=q(X)$. 
\end{prop}

\begin{proof}
Let us consider the projection $\pi_P:X\rightarrow \p^4$ 
from a general point $P$ of $\p^5$ to a 
hyperplane. $\pi_P$ is---as observed in \cite{K1}---\emph{$4$-generic} 
(see \cite{K1} for the definition).  
Therefore, we can apply Kleiman's quadruple point formula of \cite{K1},  
and from 
this we obtain formula~\eqref{eq:4tuple} in the case of a general point. 

If $P$ is not a general point, the conclusion follows by degree reasons.  
\end{proof}

Next, we pass to prove the other formula.  

\begin{prop}
Let $S\subset\p^4$ be a smooth non-scrollar surface of degree $d$, sectional 
genus $\pi$, hyperplane and canonical divisors $H$ and $K$, respectively. 
Then the number of $4$-secant lines $h$ 
of $S$ passing through 
a general point $P\in S$ is given by the formula
\begin{equation}
h=\frac{1}{6}d^3-\frac{3}{2}d^2+d(\frac{16}{3}-\pi)+4\pi+2\chi(\OO_S)-10.
\label{eq:k}
\end{equation}
\end{prop}

\begin{proof}
$h$ is actually equal to the number of triple points of the image of 
$S$ under the projection from $P$ to a hyperplane. Therefore $h$ can be 
obtained from the triple point formula for a map $f$ from a smooth surface
to $\p^3$: we simply blow-up $S$ in $P$ and then we compose the map which 
defines the blow-up $g:\bl_P(S)\rightarrow S$ with the projection, \ie 
$f:=\pi_p\circ g$. The triple point formula can be found in \cite{LEB},  
and it is, in our situation 
\begin{equation}
h=\frac{1}{6}(\tilde{d}(\tilde{d}^2-12\tilde{d}+44)+4\tilde{K}^2
-2\tilde{c_2}-3\tilde{H}\tilde{K}(\tilde{d}-8)),
\label{eq:triple} 
\end{equation}
where $\tilde{H}=g^*H-E$ is the strict transform of $H$---$E$ is the 
exceptional divisor of the blow-up---, $\tilde{K}=g^*K+E$ the canonical 
divisor of $\bl_P(S)$, $\tilde{d}=\tilde{H}^2$ and $\tilde{c_2}$ its 
topological 
Euler-Poincar\'e characteristic. 
Clearly, we have that $\tilde{d}=d-1$, $\tilde{c_2}=c_2+1$, where 
$c_2$ is the topological 
Euler-Poincar\'e characteristic of $S$,  
$\tilde{K}^2=K^2-1$ since $S$ is not a scroll, and 
$\tilde{H}\tilde{K}=2\pi-d-1$. Then, if we express the 
invariants of $S$ in terms of the basic invariants $(d,\pi,\chi(\OO_S))$,  
we get formula~\eqref{eq:k}. 
\end{proof}

\begin{rmk}
We obtained the formulae with the 
help of S. Katz and S. A. Str\o mme's Maple package ``Schubert''.
\end{rmk}


\section{Congruences of $\p^n$}\label{sec:3}

In this section we study the irreducible first order congruences  
$B$ that are given by the families of the $4$-secant lines of  smooth 
threefolds 
$X\subset\p^5$ (with the notations for its invariants given in 
Section~\ref{sec:2}), proving 
Theorem~\ref{thm:primo}. 
We need the following preliminary result: 

\begin{lem}
If $S(=X\cap H)$ is not a scroll, the following formula holds:
\begin{equation}
\begin{split}
0&=\frac{1}{8}d^4-\frac{23}{12}d^3-d^2(\pi-\frac{83}{8})
-d(\frac{355}{12}-11\pi-2\chi(\OO_S))+\\
&\quad +\frac{1}{2}\pi^2-\frac{57}{2}\pi-17\chi(\OO_S)+53. 
\label{eq:4k}
\end{split}
\end{equation}
\end{lem}

\begin{proof}
Formula \eqref{eq:4k} is formula \eqref{eq:agen}, which is, in our situation, 
\begin{equation}
4h=1+a_1 \label{eq:1}
\end{equation}
(with the notations of the proof of Theorem~\ref{thm:dgb}). 
Now, $h$ is the algebraic multiplicity of the fundamental locus on the 
variety $V_{\Pi}$ of the lines of $B$ that meet a general $3$-plane $\Pi$, 
\ie if we fix $P\in X$, there are $h$ lines of $B$ though it that meet $\Pi$ 
also. The hyperplane $\overline{P\Pi}$ intersects $X$ in a smooth surface 
$S$, and the $h$ lines are exactly the ones that are $4$-secants to $S$. 
So $h$ is given by formula~\eqref{eq:k}.

By definition, $a_1$ is the number of lines of $B$ contained in a hyperplane 
$H$ and that meet a line $\ell\subset H$, \ie it is the degree of the 
hypersurface of the $4$-secant lines of $S$. But this is formula can be 
easily deduced from \cite{Leba} and it is 
\begin{equation}
a_1=\frac{1}{8}d^4-\frac{5}{4}d^3+d^2(\frac{35}{8}-\pi)
+d(+7\pi+2\chi(\OO_S)-\frac{33}{4})+\frac{1}{2}\pi^2-\frac{25}{2}\pi
-9\chi(\OO_S)+12.\label{eq:2} 
\end{equation}

Substituting formulae \eqref{eq:k} and \eqref{eq:2}  
in \eqref{eq:1}, we get formula \eqref{eq:4k}. 
\end{proof}

\begin{rmk}
During the proof of the preceding proposition, we have calculated the first 
multidegree $a_1$ for the family (possibly reducible) of the $4$-secant lines 
of a smooth threefold $X$ of $\p^5$.  
$a_2$ is instead the number of $4$-secant lines contained in a $3$-dimensional 
linear space $G$, \ie the number of $4$-secant lines of the smooth curve 
$G\cap X$, and this formula is in \cite{LB}:
\begin{equation}
a_2=\frac{1}{12}d^4-d^3+\frac{53}{12}d^2-\frac{17}{2}d+6-\frac{1}{2}\pi d^2
+\frac{7}{2}d\pi-\frac{13}{2}\pi+\frac{1}{2}\pi^2.\label{eq:3}
\end{equation}
\end{rmk}

\begin{lem}\label{lem:lemme}
The family $B$ of the $4$-secant lines of the threefold $X$ of 
case~\textup{\eqref{case}} 
of Theorem~\textup{\ref{thm:primo}} is irreducible. 
\end{lem}

\begin{proof}
First of all, recall that $X=\cup_{\ell\in Z}\ell$, where 
$Z:=\Gr(1,5)\cap \p^8(\subset\p^{14})$: see Proposition 3.2 of 
\cite{Cha}. 

Applying Formula \eqref{eq:4tuple}, we get that $B$ is a first order 
congruence. 
If $B$ were reducible, we could write $B=B_0\cup B_1$, 
where $B_1$ is the irreducible component of $B$ of 
order $1$ and $B_0$ is the union of the other components (which have order 
zero); besides, we denote by $a_j(i)$ the $j$-th multidegree of $B_i$. 
Now, Formula~\eqref{eq:1} becomes $4h=1+a_1(1)$, but in Formula~\eqref{eq:2} 
we have to substitute $a_1$ with $a_1(0)+a_1(1)$, and therefore $a_1(0)=0$. 
So, the only possibility is $[B_0]=a_2(0)\sigma_{2,2}$, \ie 
$X$ contains at least a surface of a $\p^3$ of degree at least $4$. 

If $Y\subset\Pi\cong\p^3$ were one of these surfaces, it is enough to prove 
that 
$Y$ is not a scroll. 
In fact, once we have proven this, we can see that $Y$ would be birationally 
equivalent to $Z$, since through almost every point of $Y$ there 
is a line of $X$ not contained in $\Pi$. 
But  $Z$ is the intersection of $\Gr(1,5)$ 
with six general hyperplanes of $\p^{14}$, while in ${\p^{14}}^*$ the set 
of hyperplanes consisting of the special Schubert cycles of lines meeting 
a $3$-plane has codimension six (it is in fact $\Gr(3,5)$) and hence 
disjoint with the span of six general elements. 

If $Y$ were a scroll, then its corresponding curve $C$ in $X$ would be 
contained in a fixed $\Gr(1,\Pi)\subset\Gr(1,5)$. 
But both $X$ and $\Gr(1,\Pi)$ are linear 
sections of $\Gr(1,5)$ and therefore $C$ could be (at most) a conic, and so 
$Y$ could only be a quadric. 
\end{proof}

\begin{rmk}
In the previous proof, the fact that we can exclude the non-scrollar case, 
it has been suggested to us by the referee. Our proof was longer and based 
on the fact that the sectional genus of $Y$ would be the same of 
$Z$. 
\end{rmk}

\begin{proof}[Proof of Theorem \textup{\ref{thm:primo}}]
By Theorem \ref{thm:HC}, it is enough to consider the cases of $\p^n$ 
with $n\le 5$. 

By Theorem \ref{thm:dgb} (and by Castelnuovo's bound) we get that in 
$\p^3$ the only congruence is the one given by the secant lines of the 
twisted cubic. 

In $\p^4$, again by Theorem \ref{thm:dgb}, we obtain that the surfaces $X$ we 
are looking for have to satisfy $4\le \deg(X)\le 8$, and since the smooth 
surfaces are classified up to degree ten, see \cite{DP}, if we apply the 
triple point formula (see \cite{LEB}, or formula \eqref{eq:triple}, erasing 
the tildes) to them, we get only the cases of the theorem. 

Passing to the next case, 
since the smooth threefolds of $\p^5$ are classified  up to 
degree $12$ (see \cite{BSS}), 
we can check which of them have an apparent quadruple point, 
and it turns out that they are the ones of the list of the theorem. 
The three congruences are indeed irreducible, case~\eqref{case} by 
Lemma~\ref{lem:lemme}, the other two from their description as 
degeneracy loci. 

Next, from Theorem~\ref{thm:dgb} we have to exclude the cases of degrees 
$13$, $14$ and 
$15$. To do this, we calculate the possible invariants of these threefolds 
(see for example \cite{BSS}) 
and then we request that they have to satisfy $q(X)=1$ in 
formula~\eqref{eq:4tuple}, and equation~\eqref{eq:4k} applied to the general 
hyperplane section $S$ of $X$: since the smooth, 
non-degenerate scroll surfaces of $\p^4$ have degree at most $5$, see 
\cite{Lan}, we can apply this last formula without problems. It turns out that 
there cannot exist threefolds with these conditions.  

Finally we can calculate the multidegree for the three congruences of $\p^5$
(the cases of $\p^3$ and $\p^4$ can be easily deduced from 
Section \ref{sec:6}); $a_1$ is from formula \eqref{eq:2} and $a_2$ is 
given by formula~\eqref{eq:3}. 
\end{proof}

\begin{rmk}
We performed the calculations in the last proof with the help of a simple 
program in Maple. 
\end{rmk}

\bibliography{irr5}
\end{document}